\documentclass[12pt]{article}
\usepackage{amsmath}
\usepackage{amsfonts}
\usepackage{amssymb}
\usepackage{exscale,relsize}
\usepackage{amsthm}
\usepackage{hyperref}
\usepackage{graphicx}
\usepackage{color}
\usepackage{subfig}
\usepackage[mathscr]{eucal}

\def\bee{\begin{equation}}
\def\eee{\end{equation}}

\def\MM{\mathcal{M}}
\def \BB{{\cal B}}
\def \BBB{{\mathcal{B}_2}}

\def\DHLhksqrt#1#2{\setbox0=\hbox{$#1\sqrt{#2\,}$}\dimen0=\ht0
\advance\dimen0-0.2\ht0
\setbox2=\hbox{\vrule height\ht0 depth -\dimen0}%
{\box0\lower0.4pt\box2}}

\begin{document}

\thispagestyle{empty}
\bigskip\bigskip
\centerline{    }
\vskip 2 cm
\centerline{\Large\bf Computer experiments }
\medskip
\centerline{\Large\bf  with Mersenne  primes}
\bigskip\bigskip

\centerline{\large\sl Marek Wolf}
\bigskip
\begin{center}
Cardinal  Stefan  Wyszynski  University\\
Department of Mathematics and Sciences\\
ul. W{\'o}ycickiego 1/3, Auditorium Maximum, (room 113)\\
PL-01-938 Warsaw,   Poland\\
%{e-mail: mwolf@ift.uni.wroc.pl}
e-mail:  m.wolf@uksw.edu.pl
\end{center}

\bigskip\bigskip\bigskip
\begin{center}
{\bf Abstract}\\
\bigskip
\begin{minipage}{12.8cm}
We have calculated  on the computer the sum $\overline{\BB}_M$ of reciprocals of
all 47 known Mersenne primes with the accuracy of over 12000000 decimal digits.
Next we developed $\overline{\BB}_M$ into the continued fraction and calculated geometrical
means of the partial denominators of the continued fraction expansion of $\overline{\BB}_M$.
We get values converging to the Khinchin's constant. Next we calculated the $n$-th square
roots of the denominators of the $n$-th convergents of these continued fractions
obtaining values approaching the Khinchin-L{\`e}vy constant. These two results suggests
that the sum of reciprocals of  all Mersenne primes  is irrational, supporting
the common believe that there is an infinity of the Mersenne primes.
For comparison we have done the same procedures
with slightly modified set of 47 numbers obtaining quite different
results. Next we investigated the continued fraction whose partial
quotients are Mersenne primes and we argue that it should be transcendental.
\end{minipage}
\end{center}

\bigskip\bigskip\bigskip

%Key words: {\it Prime numbers, gaps between primes, Hardy and Littlewood conjecture }\\
% PACS numbers: {11A99, 11A41}

\vfill
\eject

\pagestyle{myheadings}

\bibliographystyle{abbrv}

\section{Introduction.}
\label{Introduction.}

The Mersenne primes $\mathcal{M}_n$ are primes of the form $2^p-1$ where $p$
must be a prime, see  e.g. \cite[Sect. 2.VII]{Ribenboim}.
%The New Mersenne Prime Conjecture states that $\mathcal{M}_n$ is prime if and only if $p$ is a prime of one of the forms $2^k \pm 1$ or $2^{2k}\pm 3$.
The set of Mersenne primes starts  with $\mathcal{M}_1=2^2-1, \mathcal{M}_2=2^3-1,
\mathcal{M}_3=2^5-1$ and only 47 primes of this form
are currently known, see Great Internet Mersenne Prime Search (GIMPS)
at www.mersenne.org.  The largest  known Mersenne prime has the
value $\MM_{47} = 2^{43112609}-1 = 3.1647026933  \ldots \times 10^{12978188}$. % 1.6987 \ldots \times 10^{12837063}$.
In general the largest known primes
are the Mersenne primes, as the Lucas--Lehmer primality test applicable only to
numbers of the form $2^p-1$ needs a multiple of $p$ steps, thus the complexity
of checking primality of $\mathcal{M}_n$ is $\mathcal{O}(\log(\mathcal{M}_n))$.
Let us remark that algorithm of Agrawal, Kayal and Saxena (AKS) for arbitrary
prime $p$  works in  about $\mathcal{O}(\log^{7.5}(p))$ steps and
modification by Lenstra and Pomerance has complexity $\mathcal{O}(\log^{6}(p))$.

\begin{figure}
\begin{center}
\includegraphics[width=0.9\textwidth, angle=0]{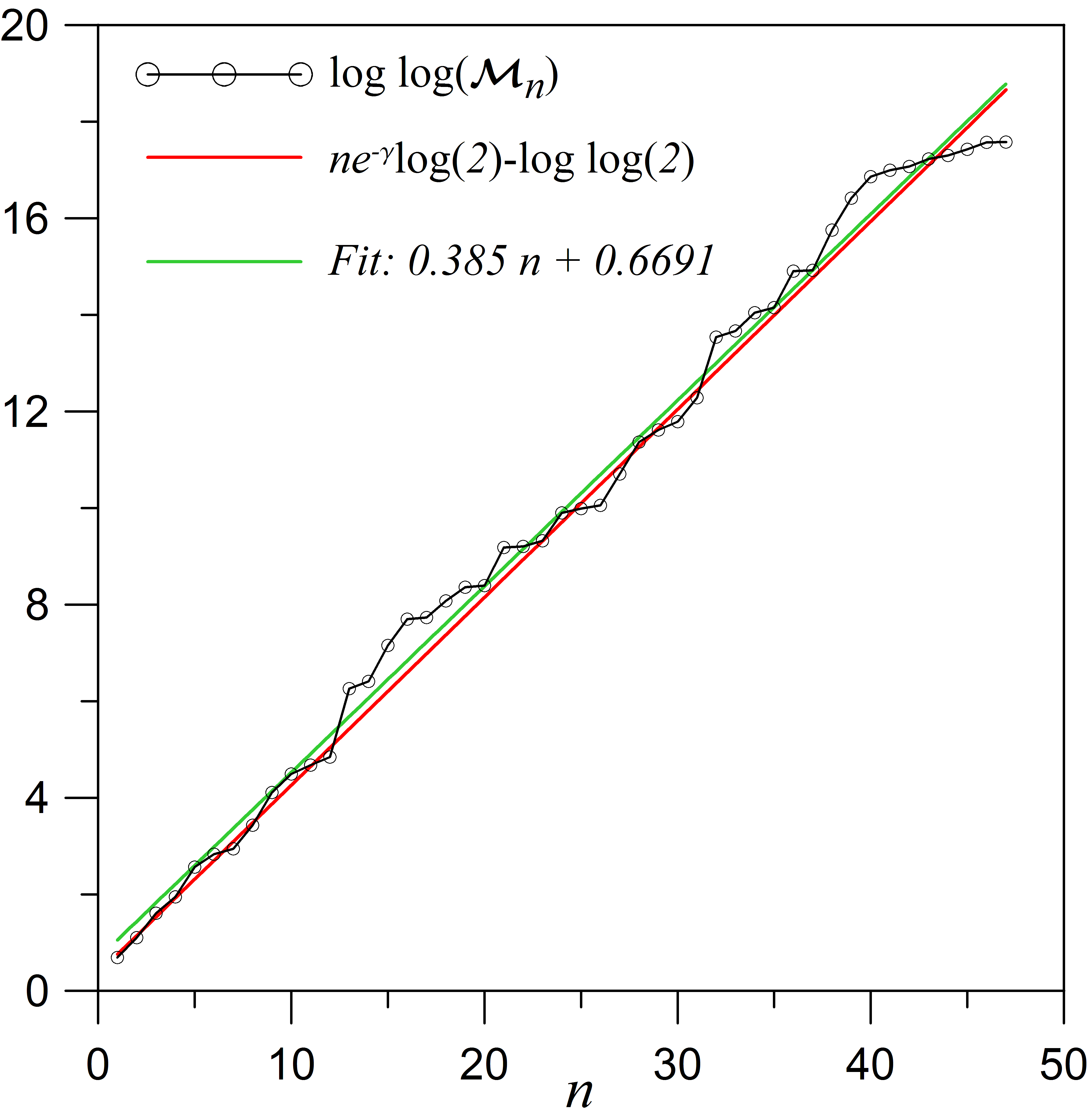} \\
\vspace{0.2cm}\caption{  The plot of $\log\log(\mathcal{M}_n)$ and the Wagstaff conjecture
(\ref{Wagstaff}).  The fit was made to all known $\mathcal{M}_n$ and it is
$0.3854n + 0.6691$, while $n e^{-\gamma}\log(2)-\log\log(2) \approx
0.3892 n +0.3665$. The rather good  coincidence of $\log\log(\mathcal{M}_n)$
and (\ref{Wagstaff}) is seeming, as to get original $\mathcal{M}_n $'s the errors
are amplified to huge values by double exponentiation.}
\label{Fig-Wagstaff}
\end{center}
\end{figure}

There is no proof of the infinitude of $\mathcal{M}_n$, but a common belief is  that as
there are presumedly infinitely many even perfect numbers thus there is also an
infinity of Mersenne primes. S. S. Wagstaff  Jr. in \cite{Wagstaff1983} (see also
\cite[\S 3.5]{Schroeder}) gave heuristic arguments, that ${\mathcal{M}_n}$ grow
doubly exponentially:
\bee
\log_2 \log_2 \mathcal{M}_n \sim n e^{-\gamma},
\label{Wagstaff}
\eee
where $\gamma=0.57721566\ldots$ is the Euler--Mascheroni constant.
In the Fig. \ref{Fig-Wagstaff} we compare the Wagstaff conjecture with all 47
presently known Mersenne primes  $\mathcal{M}_n$. Of these 47 known  $\mathcal{M}_n=2^p-1$
there are 27 with $p \mod 4=1$ and 19  with $p \mod 4=3$. It is in opposite to
the set of all primes where the phenomenon of  Chebyshev bias is known:
for initial primes there are more primes $p\equiv 3 \pmod 4$ than $p\equiv 1 \pmod 4$,
\cite{Kaczorowski}, \cite{Rubinstein_Sarnak}.

In this paper we are going to exploit two facts about the continued fractions
to support the conjecture on the infinitude of Mersenne primes: the
existence of the   Khinchin constant and  Khinchin–--L{\`e}vy constant. We calculate
the sum of reciprocals of the Mersenne primes $\BB_M=\sum_n 1/\MM_n$; if there is infinity
of Mersenne primes then this number $\BB_M$ should be irrational (at least, because
it is probably even transcendental, as it is difficult to imagine the polynomial
with some mysterious integer coefficients whose one of roots should be $\BB_M$).

There exists a method based on the continued fraction expansion which
allows  to detect whether a given number $r$ can be irrational or not. Let
\bee
r = [a_0(r); a_1(r), a_2(r), a_3(r), \ldots]=
a_0(r)+\cfrac{1}{a_1(r) + \cfrac{1}{a_2(r) + \cfrac{1}{a_3(r)+\ddots}}}
\eee
be the continued fraction expansion of the real number $r$, where $a_0(r)$ is an
integer and all $a_k(r)$ with $k\geq 1$  are positive integers.
The quantities $a_k(r) $ are called partial quotients or  the partial denominators.
Khinchin has proved \cite{Khinchin}, see also \cite{Ryll-Nardzewski1951}, that
\bee
\lim_{n\rightarrow \infty} \big(a_1(r) \ldots a_n(r)\big)^{\frac{1}{n}}=
\prod_{m=1}^\infty {\left\{ 1+{1\over m(m+2)}\right\}}^{\log_2 m} \equiv K \approx 2.685452001
\label{Khinchin}
\eee
is a constant for almost all real $r$ \cite[\S 1.8]{Finch} (the term  $a_0$ is skipped
in (\ref{Khinchin})).  The exceptions are of the
Lebesgue measure zero and include {\it rational numbers}, quadratic irrationals and
some irrational numbers too, like for example the
Euler constant $e=\lim_{n\to \infty} (1+\frac{1}{n})^n=2.7182818285\ldots$ for
which the $n$-th geometrical mean tends to
infinity like $\sqrt[3]{n}$, see \cite[\S 14.3 (p.160)]{Havil03}.
The constant $K$ is called the Khinchin constant. If the sequence % quantities
\bee
K(r; n)=\big(a_1(r) a_2(r) \ldots a_n(r)\big)^{\frac{1}{n}}
\label{Kny}
\eee
for a given number $r$ tend to $K$ for $n\to \infty$ we can regard it as an indication
that $r$ is irrational  --- all rational numbers have finite number of partial
quotients in the  continued fraction expansion and hence starting with some $n_0$
for all $n>n_0$ will be $a_n=0$. It seems to be possible to construct  a sequence
of rational numbers such that the geometrical means of partial quotients of their
continued fraction  will  tend to the Khinchin constant.

The  Khinchin---L{\`e}vy's constant arises in the following way:  Let the rational
$P_n(r)/Q_n(r)$ be the $n$-th  partial convergent of the continued fraction of $r$:
\bee
\frac{P_n(r)}{Q_n(r)}=[a_0(r); a_1(r), a_2(r), a_3(r), \ldots, a_n(r)].
\label{eq-convergent}
\eee
In 1935 Khinchin \cite{Khinchin-1935} has proved that for almost all real numbers $r$
the denominators of the  finite continued fraction approximations fulfill:
\bee
\lim_{n \rightarrow \infty} \big(Q_n(r)\big)^{1/n} \equiv \lim_{n \rightarrow \infty} (L(r; n) =  L
% = e^{\pi^2/12\ln2} \equiv L = 3.275822918721811\ldots
\eee
and in  1936 Paul Levy
\cite{Levy1936} found an explicit expression for this constant $L$:
\bee
\lim_{n \to \infty}\sqrt[n]{Q_n(r)}=e^{\pi^2/12\log(2)}\equiv L=3.27582291872\ldots
\label{Levy}
\eee
$L$ is called the  Khinchin---L{\`e}vy's constant  \cite[\S 1.8]{Finch}.
Again the set of exceptions to the above limit is  of the Lebesgue measure zero
and it includes rational numbers, quadratic irrational etc.

\section{First experiment}

Let us define the sum of reciprocals of all Mersenne primes:
\bee
\BB_M=\sum_{n=1}^\infty  \frac{1}{\MM_n},
\label{def-B_M}
\eee
which  can regarded  as the analog of the Brun's constant, i.e. the sum of
reciprocals of all twin primes:
\bee
\BBB=\left( {1 \over 3} + {1\over 5}\right) + \left( {1 \over 5} + {1\over
7}\right) + \left( {1 \over 11} + {1\over 13}\right) + \ldots.
\label{def-B2}
\eee
In  1919 Brun \cite{Brun} has shown that this constant $\BBB$ is finite,
thus leaving the problem of infinity of twin primes not decided.
Today's best  numerical value is $\BBB\approx 1.90216058$, see \cite{Nicely_Brun},
\cite{Sebah-Brun}. Yet it is possible to prove that there is infinity of twins
by showing that Brun's constant is irrational \cite{Wolf-Brun} (we believe it
is even transcendental).

Using PARI \cite{PARI}  we have calculated the sum of reciprocals of all known 47
Mersenne primes $\overline{\BB}_M$ with accuracy over 12 millions digits:
\bee
\overline{\BB}_M = 0.5164541789407885653304873429715228588159685534154197\ldots  .
\eee
This  number is not recognized by the  Symbolic Inverse Calculator
(http://pi.lacim.uqam.ca) maintained by Simone Plouffe.
The bar over $\BB_M$ denotes the finite (at present consisting of
47 terms) approximation to the full sum defined in
(\ref{def-B_M}). It is not known, whether there are  Mersenne prime numbers with
exponent $20996011<p<43112609$ --- currently confirmed by GIMPS is that
$2^{20996011}-1$ is the 40-th Mersenne prime $\MM_{40}$--- it is not known whether any
undiscovered Mersenne primes exist between the 40th $\MM_{40}$ and
the 47th Mersenne prime $\MM_{47}$.  We have taken 12000035 digits of
$\overline{\BB}_M$ --- it means that we assume that there are no unknown Mersenne primes
with $p<39863137$. Using the incredibly fast procedure {\tt ContinuedFraction[$\cdot$]}
implemented in Mathematica$^\copyright$ we  calculated
the continued fraction expansion of $\overline{\BB}_M$.  The result was built from
11645012 partial denominators  $a_1=1, a_2=1, a_3=14, \ldots, a_{11645012}=4$.
The $n$-th convergent  ${P_n(r)/Q_n(r)}$, see (\ref{eq-convergent}),
approximate the value of $r$  with accuracy at least $1/Q_k Q_{k+1}$
\cite[Theorem 9, p.9]{Khinchin}:
\bee
\left|r - \frac{P_k}{Q_k} \right |<\frac{1}{Q_k Q_{k+1}}<\frac{1}{Q_k^2a_{k+1}}<\frac{1}{Q_k^2}.
\label{error}
\eee
From this it follows that if  $r$ is known with precision of  $d$ decimal  digits we can
continue with calculation of quotients $a_n$ up to such $n$ that
the denominator of the n-th convergent $Q_n^2<10^d$. We have checked that
$Q_{11645012}= 4.291385\times 10^{6000016}$.

The largest denominator was  $a_{9965536}=716699617$. We have checked correctness of the
continued fraction expansion of $\overline{\BB}_M$ by calculating backwards from
$[0; 1, 1, 14, \ldots, 4]$ the partial convergent. The Mathematica$^\copyright$ has
build in the procedure {\tt FromContinuedFraction[$\cdot$]},  but we have used our
own procedure written in PARI and implementing the recurrence:
\bee
P_{n+1} =a_nP_{n}+P_{n-1},  ~~~~~~~~ Q_{n+1} =a_nQ_{n}+Q_{n-1}, ~~~~ n\geq 1
\eee
with initial values
\bee
P_0=a_0,~~~~Q_0=1,~~~~P_1=a_0a_1+1, ~~~~Q_1=a_1.  %  ~~~~P_2=a_0a_1a_2+a_2+a_0,~~~~Q_2=a_1a_2+1
\eee

We have obtained the ratio of two mutually prime 6000018 decimal digits long integers
(it means denominator was of the order $10^{6000018}$ and hence its square was smaller
than $10^{12000035}$,  see eq.(\ref{error})):
\[
\begin{split}
~~~~~~~~~~~6000018~~digits~~~~~~~~~~~~~~\\
\frac{\overbrace{\color{white}{\frac{  }{\color{black}{2216304109121123313251143869\ldots 2210}}}}}
                {4291385759849224534616716035\ldots 2813}\\
     \end{split}
\]
whose ratio had 12000033 digits the same as $\overline{\BB}_M$.   The decimal
expansion of $\overline{\BB}_M=P/Q$ is of course periodic (recurring), see
\cite[Th. 135]{H-W},
but the length of the period is much larger than $1.2\times 10^{12}$. According
to the Theorem 135 from \cite{H-W} the period $r$ of the decimal expansion of
$\overline{\BB}_M$ is equal to the order of 10 mod $Q$, i.e. it is the smallest %determined from the equation:
positive $r$ for which
\bee
10^r \equiv 1 (\!\!\!\!\! \mod Q).
\eee
Because $Q$ being the product of all 47 Mersenne primes is of the order
$3.509\ldots \times 10^{86789810}$, we expect that the value of $r$ is much larger
than $10^{12}$. Another argument
is that  we got over 11 500 000 partial quotients of the continued fraction
of  $\overline{\BB}_M$ --- the numbers with periodic
decimal expansions have only finite number  of partial quotients different from zero.

\begin{figure}[h]
\begin{center}
\includegraphics[width=0.9\textwidth, angle=0]{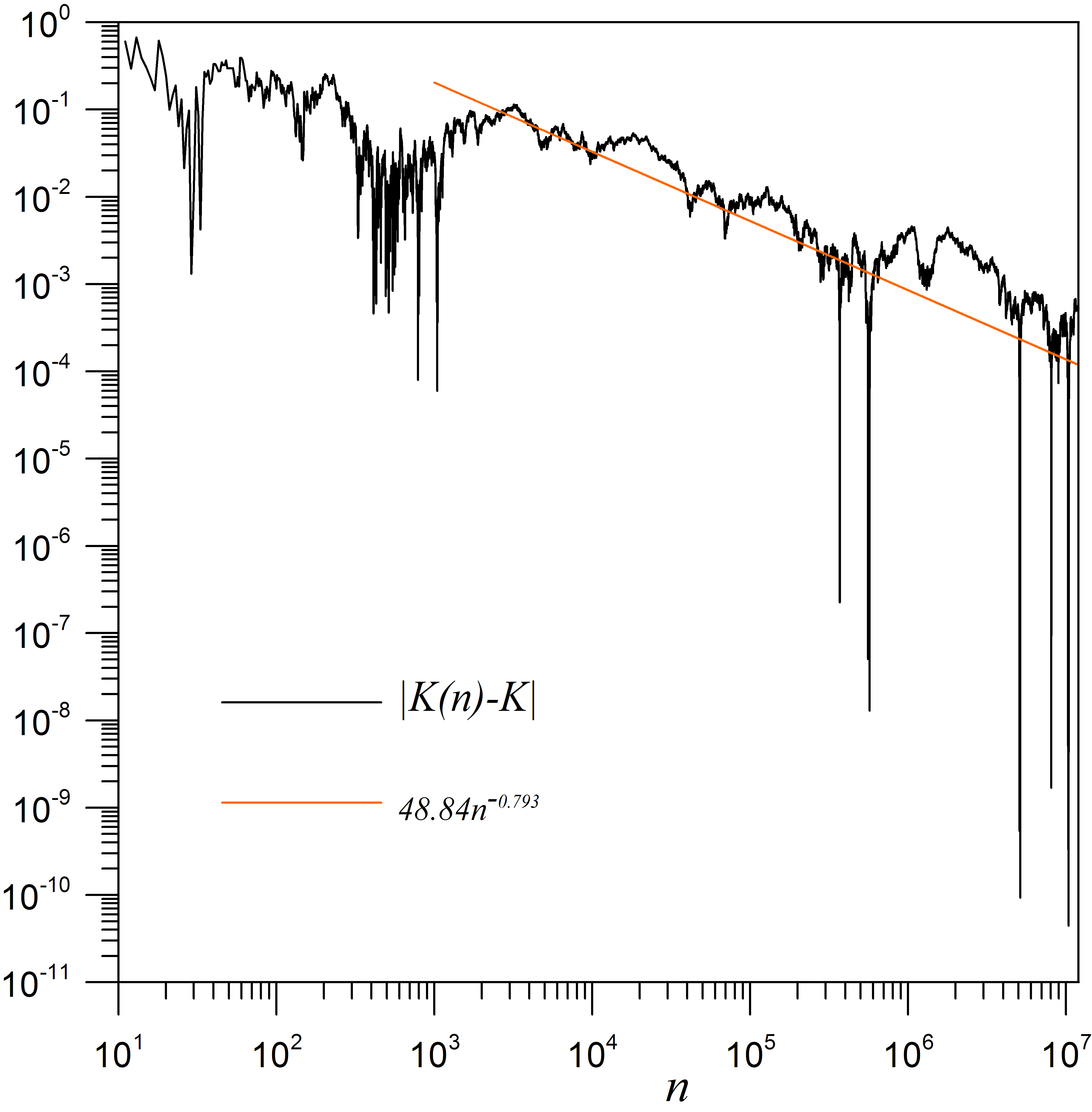}
\vspace{0.2cm}
\caption{The plot showing the distance to $K$ of the running geometrical averages
$K(n)=(a_1 a_2 \cdots a_n)^{1/n}$ for the continued fraction of $\BB_M$.  }
\label{fig-errors-K}
\end{center}
\end{figure}

From the sequence of partial quotients $a_1=1, a_2=1, a_3=14, \ldots, a_{11645012}=4$
we have calculated running geometrical means
\bee
K(n)= \left( \prod_{k=1}^{n} a_k \right)^{1/n}
\eee
for $n=11, \ldots, 11645012$. The obtained numbers $K(n)$ quickly tend to the
Khinchine constant thus in Fig. \ref{fig-errors-K} we have plotted the differences
$|K(n)-K|$. The power fit to the values for $n=1000\ldots 11645012$ gives the
decrease of the form $|K(n)-K| \sim n^{-0.79}$ and it suggests that indeed
$\lim_{n\to \infty} K(n)=K$ and thus $\BB_M$ is irrational. Indices $n$ for which the
geometric means $K(n)$ produce progressively better approximations to Khinchin's
constant are:
\bee
1, 3, 2, 16, 17, 21, 24, 26, 29, 412, 788, 1045, 369625, 369636, \ldots ,
5137093, 10389989;
\eee
the smallest value of $|K(n)-K|$ was $4.455957\ldots\times 10^{-11}$. This sequence
can be regarded as the counterpart  to the A048613 at OEIS.org.

\begin{figure}[h]
\begin{center}
\includegraphics[width=0.9\textwidth,angle=0]{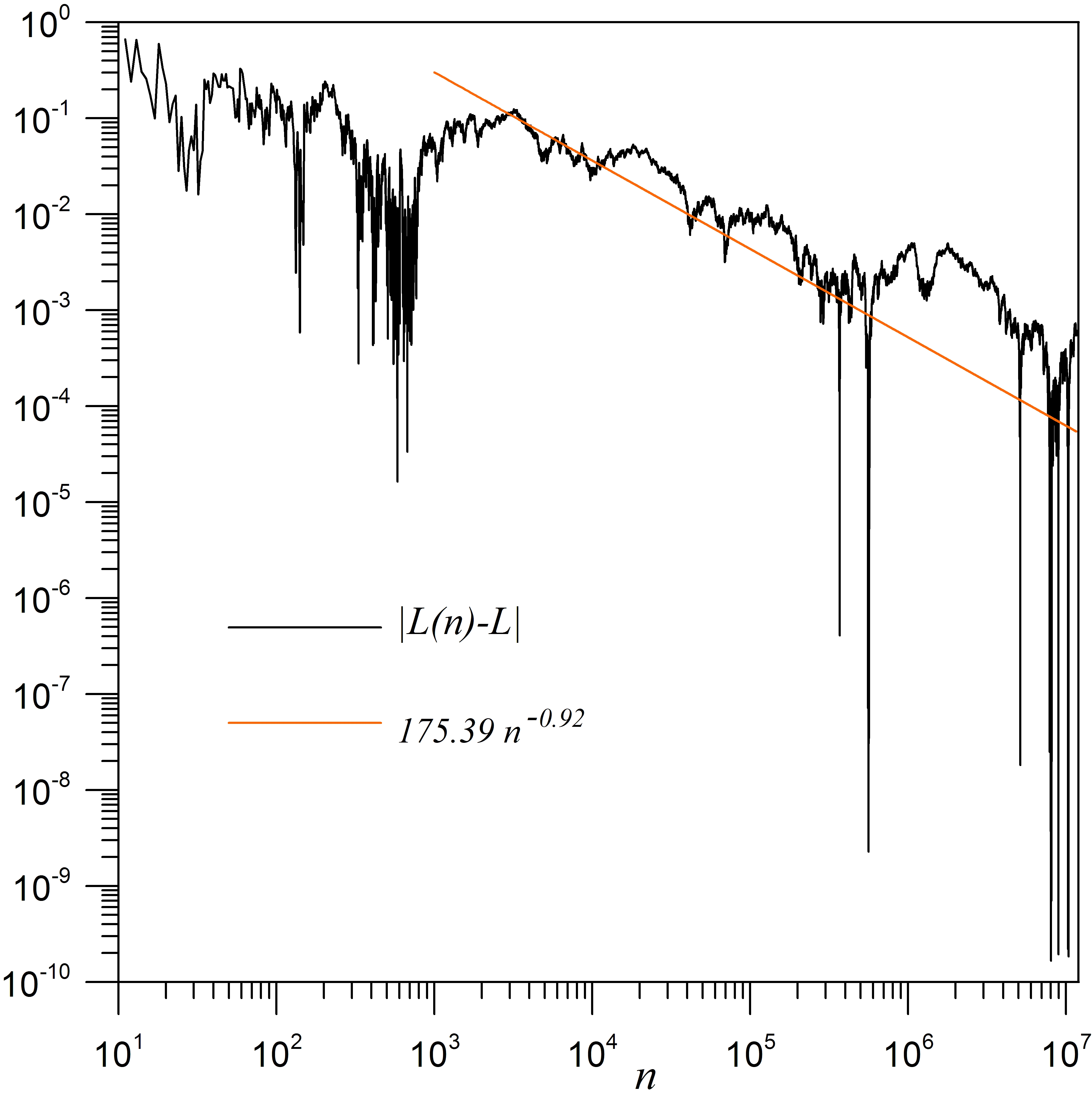}
\vspace{0.2cm}
\caption{The plot showing the distance to $L$ of the $(Q(n))^{1/n}$ obtained
from the  partial convergents of the  continued fraction of $\BB_M$
for $n=11, \ldots, 11645012$. }
\label{fig-errors-L}
\end{center}
\end{figure}

Next we calculated running (i.e. for $n=11, \ldots, 11645012$)  partial quotients
$P_n/Q_n$ and then the
quantities $L(n)=\sqrt[n]{Q_n}$, which for almost all irrational numbers should tend to
the Khinchine--Levy constant.  The behaviour of $\sqrt[n]{Q_n}$ is  shown in
Fig.\ref{fig-errors-L}. Again we see that these quantities  tend to the
limit $L$; the fitting of the power--like dependence for $n>10$ gives that $|L(n)-L|\approx
175.39 n^{-0.92}$.  The shape of the plot in this figure is similar to the plot of
$|K(n)-K|$  in Fig. \ref{fig-errors-K}.

\begin{figure}
\vspace{-1.0cm}
\hspace{-3.5cm}
\begin{center}
\includegraphics[height=0.4\textheight,angle=0]{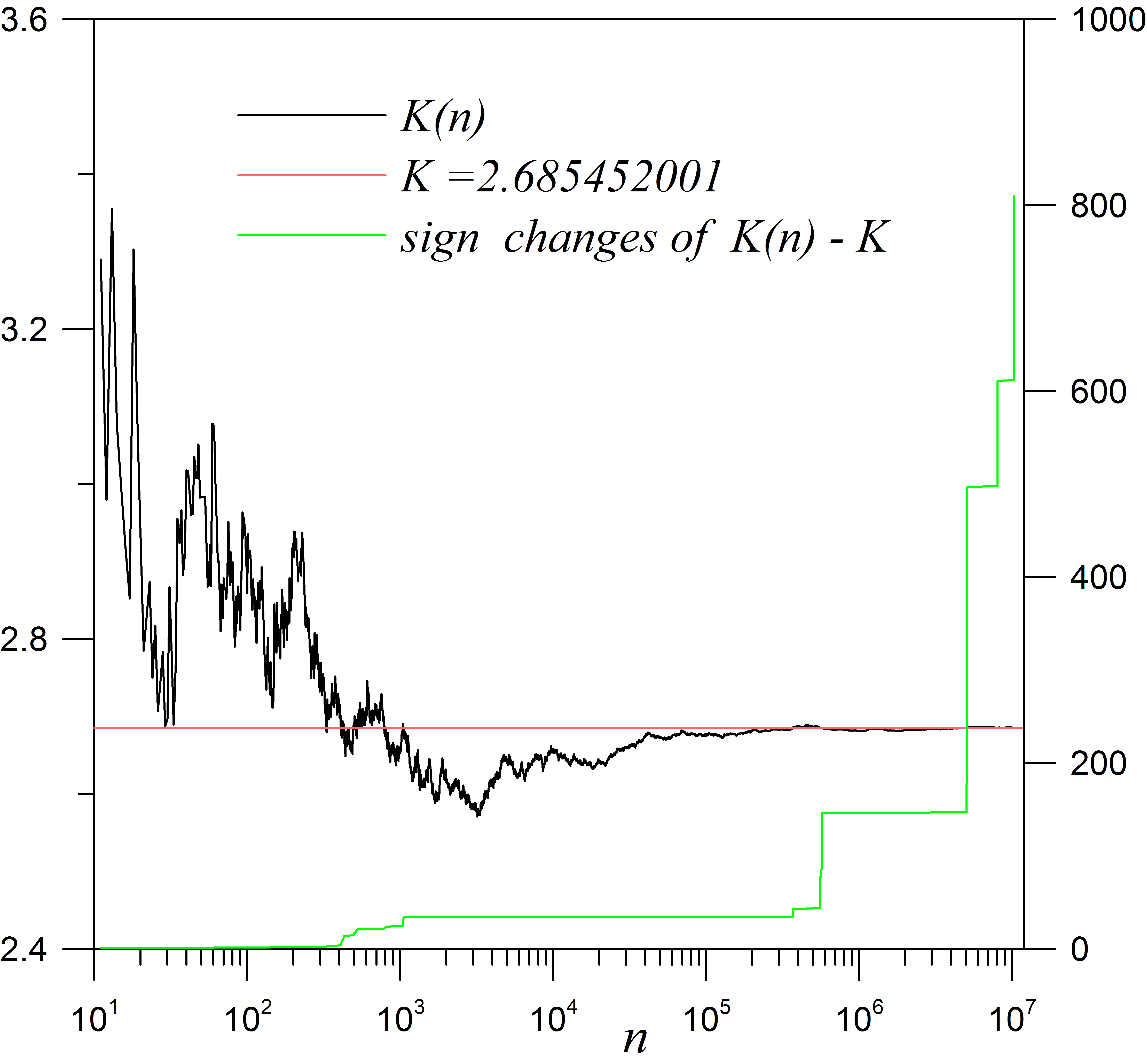} \\
\vspace{0.2cm}
\caption{The plot of $K(n)$ in black approaching the Khinchine constant $K=2.685452\ldots$
(in red) with values presented on left $y$-axis.  In green are presented numbers of
sign changes of $K(n)-K$ up to $n$ --- the right $y$-axis is for this plot.}
\label{fig-Mersenne_Khinchine_zmiany_znaku}
\end{center}
\end{figure}

\begin{figure}
\vspace{-1.3cm}
\hspace{-3.5cm}
\begin{center}
\includegraphics[height=0.4\textheight,angle=0]{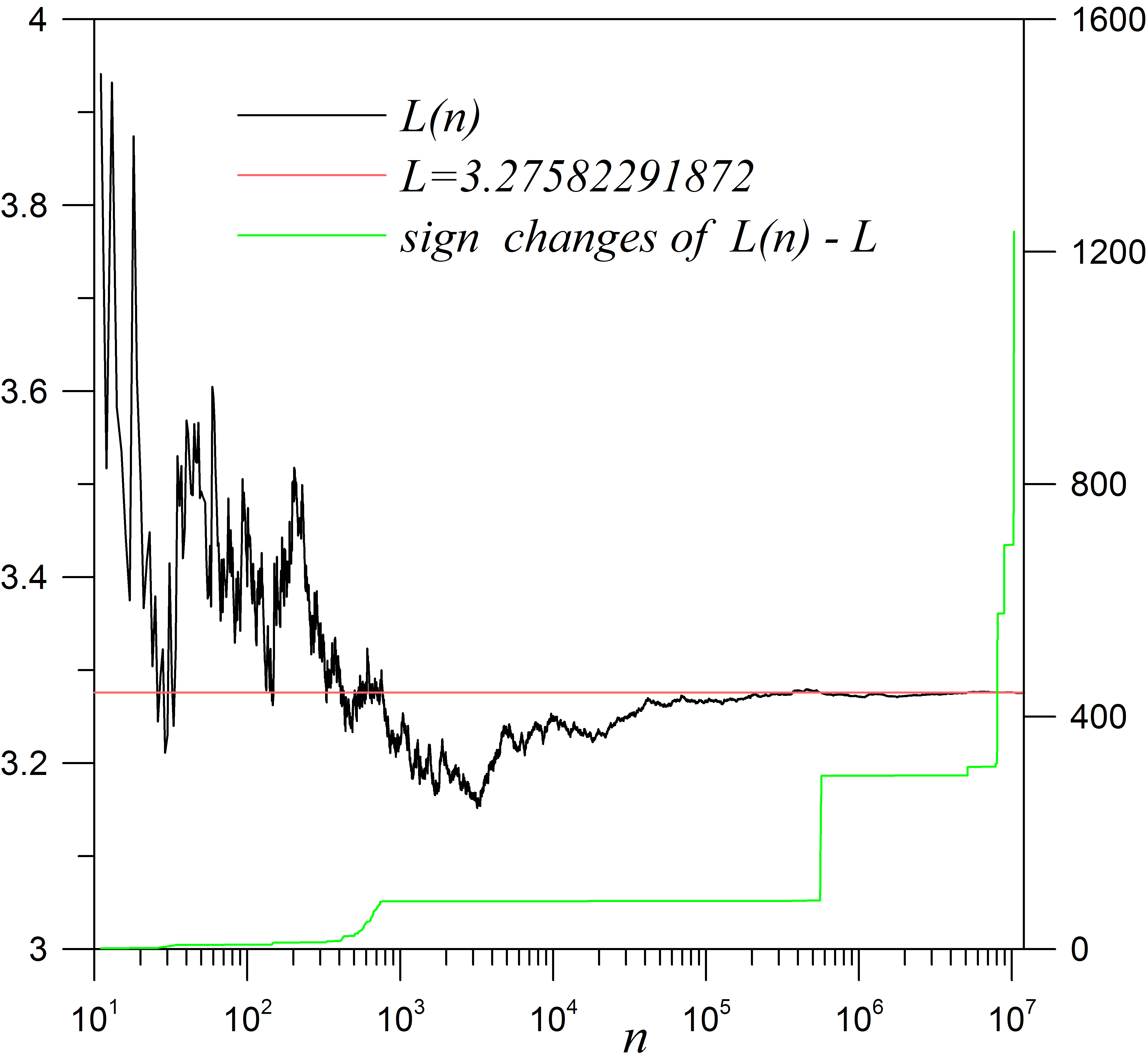} \\
\vspace{0.2cm}
\caption{ The plot of $L(n)$ in black approaching the Khinchine--Levy constant
$L=3.27582291872\ldots$ (in red) with values presented on left $y$-axis.
In green are presented numbers of  sign changes of $L(n)-L$ up to $n$ ---
the right $y$-axis is for this plot.}
\label{fig-Mersenne_Levy_zmiany_znaku}
\end{center}
\end{figure}

Both differences $K(n)-K$ and $L(n)-L$ have a lot sign changes for $n<11645012$.
Figures \ref{fig-Mersenne_Khinchine_zmiany_znaku} and \ref{fig-Mersenne_Levy_zmiany_znaku}
present the plots of these differences together with the number of sign changes.

The data presented in Figures  \ref{fig-errors-K} and  \ref{fig-errors-L}
provide the hints that $\BB_M$ is irrational and hence that there is infinity
of Mersenne primes.   But we are convinced
$\BB_M$ {\it is in fact transcendental}. In favor of this claim we recall here the result
of A. J. van der Poorten and J. Shallit \cite{PoortenShallit} that the following sum
\bee
\frac{1}{2^1}+\frac{1}{2^2}+\frac{1}{2^3}+\frac{1}{2^5}+\ldots+\frac{1}{2^{F_n}}+\ldots
\eee
where $F_n$ are Fibonacci numbers, is transcendental. %Let us remark,
It is well known that the Liouville number
\bee
\frac{1}{2^{1!}}+\frac{1}{2^{2!}}+\frac{1}{2^{3!}}+\frac{1}{2^{4!}}+\dots+\frac{1}{2^{n!}}+\ldots
\eee
is transcendental see \cite[Theorem 192]{H-W}.  In $\BB_M$, assuming the
Wagstaff conjecture,  unfortunately
the terms  decrease slower: $n!>2^{n}>2^{e^{-\gamma}n}$ for $n\geq 4$ but faster than
$F_n=\bigg\lfloor\frac{\varphi^n}{\sqrt 5} + \frac{1}{2}\bigg\rfloor$,
where  $\varphi = \frac{1 + \sqrt{5}}{2} \approx 1.6180339887\dots $.

\begin{figure}[h]
\begin{center}
\includegraphics[width=0.9\textwidth,angle=0]{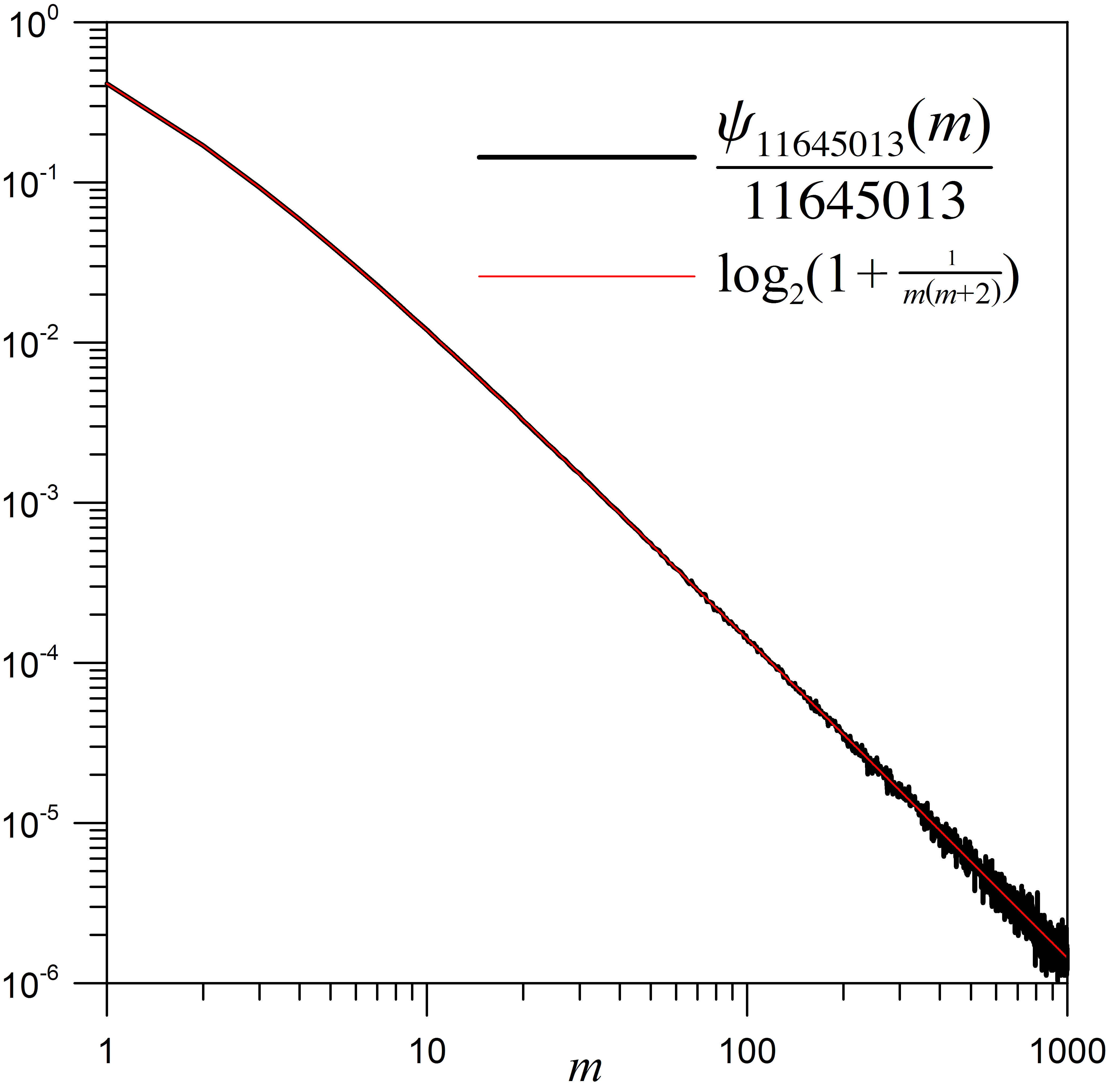}
\vspace{0.2cm}
\caption{The plot of the measured for the continued fraction of $\BB_M$ probability
to find the partial quotient  $a_k=m$ for the continued fraction of $\BB_M$. }
\label{fig-Kuzmin}
\end{center}
\end{figure}

Let $\psi_n(m)$ denotes the number of partial quotients $a_k$ with $k=1, 2, \ldots, n$
which are equal to $m$:
\[
\psi_n(m)=\sharp\{ k:  k\leq n \mbox{ and }a_k=m\}.
\]
Then the Gauss--Kuzmin theorem (for excellent exposition see e.g.
\cite[\S 14.3]{Havil03})  asserts that
\bee
\lim_{n \to \infty} \frac{\psi_n(m)}{n}=\frac{\log\left(1+\frac{1}{m(m+2)}\right)}{\log(2)}
\eee
for continued fractions of almost all real numbers. In other words, the probability
to find the partial quotient $a_k=m$ is equal to $\log_2(1+1/m(m+2))$. In Fig.
\ref{fig-Kuzmin} we present the plot of the $\frac{\psi_{11645013}(m)}{11645013}$
for the continued fraction of $\overline{\BB}_M$ and $m=1, 2, \ldots 1000$ together with
prediction given by  the Gauss--Kuzmin theorem finding excellent agreement.

Finally let us notice, that the number $\overline{\BB}_M$ computed with 12000035 digits
is normal in the base 10, see  Table I.

\vskip 0.4cm
\begin{center}
{\sf TABLE {\bf I}}\\
Illustration of  the normality of $\overline{\BB}_M$: the numbers in second row
gives the number of digits $ 0, 1, \ldots 9$ appearing in the decimal expansion of
$\overline{\BB}_M$  and the third row contains the ratio of numbers in second row
divided by 12000035.
\\ \bigskip
\begin{tabular}{|c|c|c|c|c|} \hline
 0  &  1 &  2 &  3 &  4   \\  \hline
1200553   & 1199322 & 1199420 & 1200548 & 1199397   \\ \hline
0.1000458   & 0.0999432 & 0.0999514 & 0.1000454 & 0.0999495   \\ \hline
\hline
 5  &  6 &  7 &  8 &  9   \\  \hline
1198596   & 1200876 & 1200056 & 1201757 & 1199510   \\ \hline
0.0998827   & 0.1000727 & 0.1000044 & 0.1001461 & 0.0999589   \\ \hline
\end{tabular}\\
\end{center}

\begin{figure}
\vspace{-2.3cm}
\hspace{-3.5cm}
\begin{center}
\includegraphics[height=0.4\textheight,angle=0]{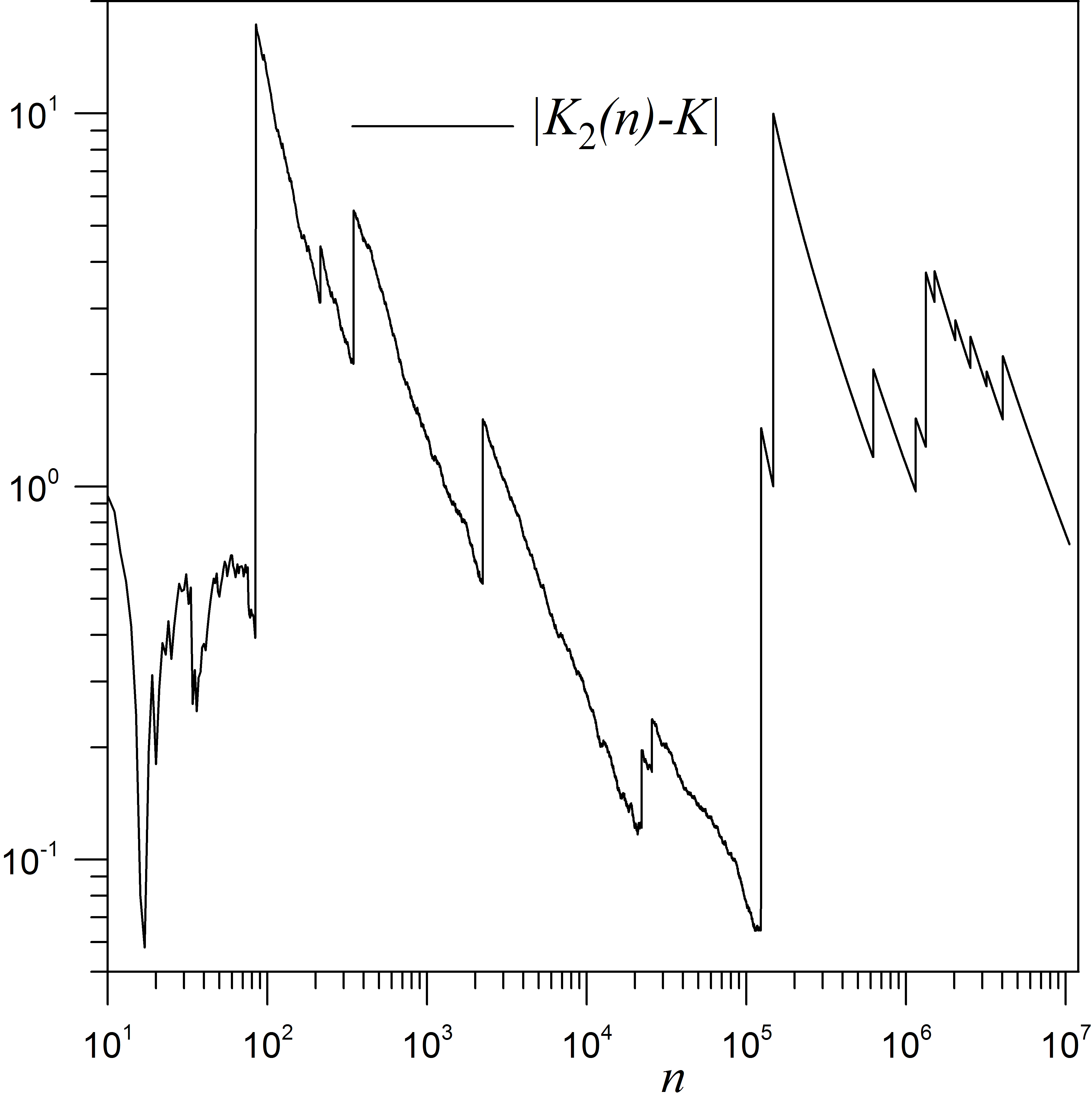} \\
\vspace{0.2cm}
\caption{ The plot showing the distance to $K$ of the running geometrical averages
$K_2(n)$ for the continued fraction of $\mathcal{S}$
for $n=11, \ldots, 10550114$.}
\label{fig-sztuczny_2_CF_Khinchin}
\end{center}
\end{figure}

\begin{figure}
\vspace{-2.3cm}
\hspace{-3.5cm}
\begin{center}
\includegraphics[height=0.4\textheight,angle=0]{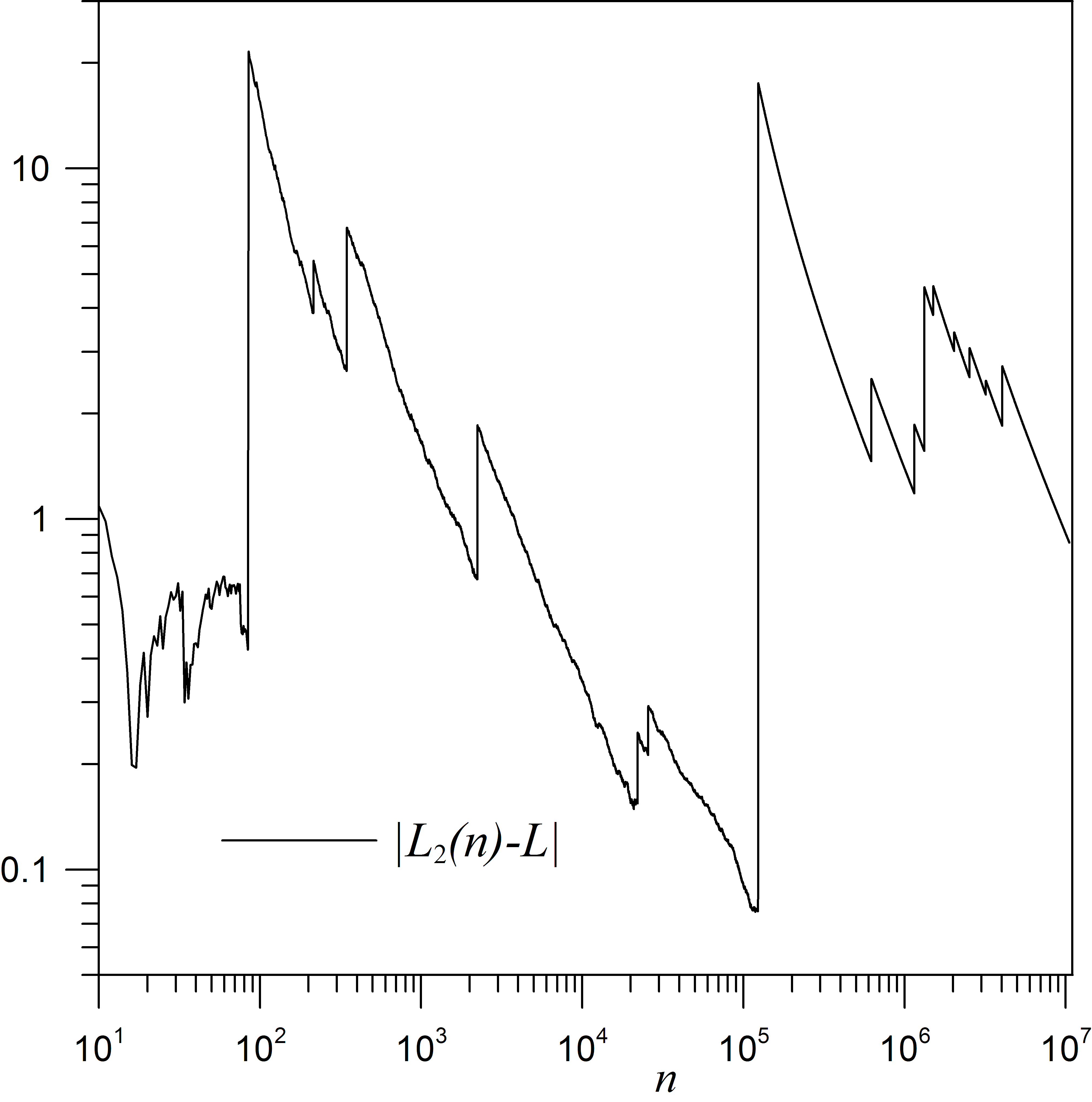} \\
\vspace{0.2cm}
\caption{ The plot showing the distance to $L$ of the $(Q_2(n))^{1/n}$ obtained
from the  partial convergents of the  continued fraction of $\mathcal{S}$
for $n=11, \ldots, 10550114$. }
\label{fig-sztuczny_2_CF_Levy}
\end{center}
\end{figure}

For comparison we have repeated the above procedure for  artificial  set
of 47 numbers of the size corresponding to known  Mersenne primes. We have simply
skipped -1 in the  Mersenne primes and using PARI we have computed with over
120000000 digits the sum:
\[
\mathcal{S}=\frac{1}{2^2} +\frac{1}{2^3} +\ldots +\frac{1}{2^{42643801}} +\frac{1}{ 2^{43112609}}
\]
This number $\mathcal{S}$ is the ratio of the form $A/2^{43112609}$, where
$\gcd(A, 2{43112609})=1$.
From \cite[\S 9.2]{H-W} we know that $\mathcal{S}$ has {\it  terminating} decimal expansion
consisting of 43112609 decimal digits, thus calculating 12000000 digits of this sum
makes sens as it does not contain recurring periodic patterns of digits. We have
developed $\mathcal{S}$ into the continued fraction, what resulted in %  The obtained
10550114 partial quotients. The calculated quantities for this case we denote
with the subscript 2: $Q_2(n), K_2(n), L_2(n)$ to distinguish them  from earlier
experiment for $\overline{\BB}_M$.  We have  calculated
the running geometrical averages of the partial quotients $K_2(n)$ and the results
are presented in Figure \ref{fig-sztuczny_2_CF_Khinchin}.
Next we calculated 10550114 partial convergents $P_2(n)/Q_2(n),  n=1, 2, \ldots,
10550114$ and from them the quantities  $L_2(n)\equiv (Q_2(n))^{1/n}$, which should
tend to the Levy constant $L$. In  Figure \ref{fig-sztuczny_2_CF_Levy} the
differences $|L_2(n)-L|$ are plotted. Obtained plots are completely different
from those seen in  Figures  \ref{fig-errors-K} and  \ref{fig-errors-L}
and they suggest  $K_2(n)$ as well as $L_2(n)$ do not have the limit.
In this artificial case we have encountered the phenomenon of extremely
large partial denominators:  there were $a_n$ of the order  $10^{70548}, 10^{97732}$
and $10^{279910}$. These large partial denominators are responsible for the smaller
number of $a_k$ than for $\overline{\BB}_M$, see (\ref{error}).

\section{Second experiment}

Let us define the supposedly infinite and convergent continued fraction
$u_{\mathcal{M}}$ by taking $a_n=\mathcal{M}_n$:
\bee
u_{\mathcal{M}}=[0; 3,~ 7,~ 31,~ 127,~ 8191,~ 131071,~ 524287,~ 2147483647,~ \ldots]
\label{def-uM}
\eee

Using all 47 Mersenne primes  $3, 7, 31, \ldots, 2^{43112609}-1$
in a couple of minutes we have calculated $u_{\mathcal{M}}$
with the precision of 10000000 digits;  first 50 digits of $u_{\mathcal{M}}$ are:
\bee
u_{\mathcal{M}}=0.31824815840584486942596202748140694243806236564\ldots
\eee
This  number is not recognized by the  Symbolic Inverse Calculator
(http://pi.lacim.uqam.ca).
Because $1/Q_{47}^2(u_{\mathcal{M}})\approx 1.84313\times 10^{-173579621}$   %  0.542555\times 10^{173579620}$
it follows from (\ref{error}) that theoretically it is possible to
obtain the value of $u_{\mathcal{M}}$ from presently known 47
Mersenne primes with over 170,000,000 decimal  digits of accuracy.
Of course $u_{\mathcal{M}}$ is  the exception to the Khinchin and Levy Theorems
in view of the very fast growth of $u_{\mathcal{M}}$ --- see the Wagstaff \cite{Wagstaff1983}
conjecture (\ref{Wagstaff}).

There is a vast literature concerning the transcendentality of continued fractions.
For example the continued fraction
\bee
[0; 2^{1!}, 2^{2!}, 2^{3!}, \ldots, 2^{n!},\ldots]
\eee
is transcendental, see \cite[Theorem 192]{H-W}, \cite[Example 4]{Sondow-2004}.

The Theorem of H. Davenport and K.F. Roth \cite{Davenport1955} states, that if
the denominators $Q_n$ of convergents of the continued fraction
$r=[a_0; a_1, a_2, \ldots]$ fulfill
\bee
\limsup_n \frac{\sqrt{\log(n)} \log (\log (Q_n(r)))}{n}=\infty
\label{Davenport-Roth}
\eee
then $r$ is transcendental. This  theorem requires for the transcendence
of $r$ very fast increase of denominators of the convergents: at least doubly
exponential growth
%something like \bee
%Q_n>e^{e^{nf(n)}}, ~~~~f(n) \to \infty    \eee
is required for (\ref{Davenport-Roth}). The set of continued fractions which can satisfy
the  Theorem of H. Davenport and K.F. Roth is of measure zero, as it follows from the
Theorem 31 from the Khinchin's book \cite{Khinchin}, which asserts there exists an absolute
constant $B$ such that for {\it almost all} real numbers $r$ and sufficiently large
$n$ the denominators of its  continued fractions  satisfy:
\bee
Q_n(r)<e^{Bn}.
\eee

The paper of A. Baker \cite{Baker1962} from 1962 contains a few theorems
on the  transcendentality of  Maillet type  continued fractions \cite{Maillet},
i. e.  continued fractions with bounded partial quotients which have transcencendental
values. Besides Maillet  continued fractions there are some specific families of
other continued fractions of which it is known that they are
transcendental. In the papers \cite{Queffelec}, \cite{Adamczewski2007-AMM} it was
proved that the Thue--Morse
continued fractions with bounded partial  quotients are transcendental.
Quite recently there appeared the preprint
\cite{Bugeaud2010} where the transcendence of the Rosen continued fractions was
established. For more examples see \cite{Adamczewski-Bugeaud-Davison}.

In the paper \cite{Adamczewski2007} B.  Adamczewski and  Y. Bugeaud,
among others, have improved (\ref{Davenport-Roth}) to the form:
\bee
{~~~~\rm If}~~~\limsup_n \frac{ \log (\log (Q_n(r)))}{n^{2/3}\log(n)^{2/3}\log(\log(n))}=\infty
\label{Adamczewski-Bugeaud}
\eee
then $r$ is transcendental.

\begin{figure}
\vspace{-2.3cm}
\hspace{-3.5cm}
\begin{center}
\includegraphics[height=0.4\textheight,angle=0]{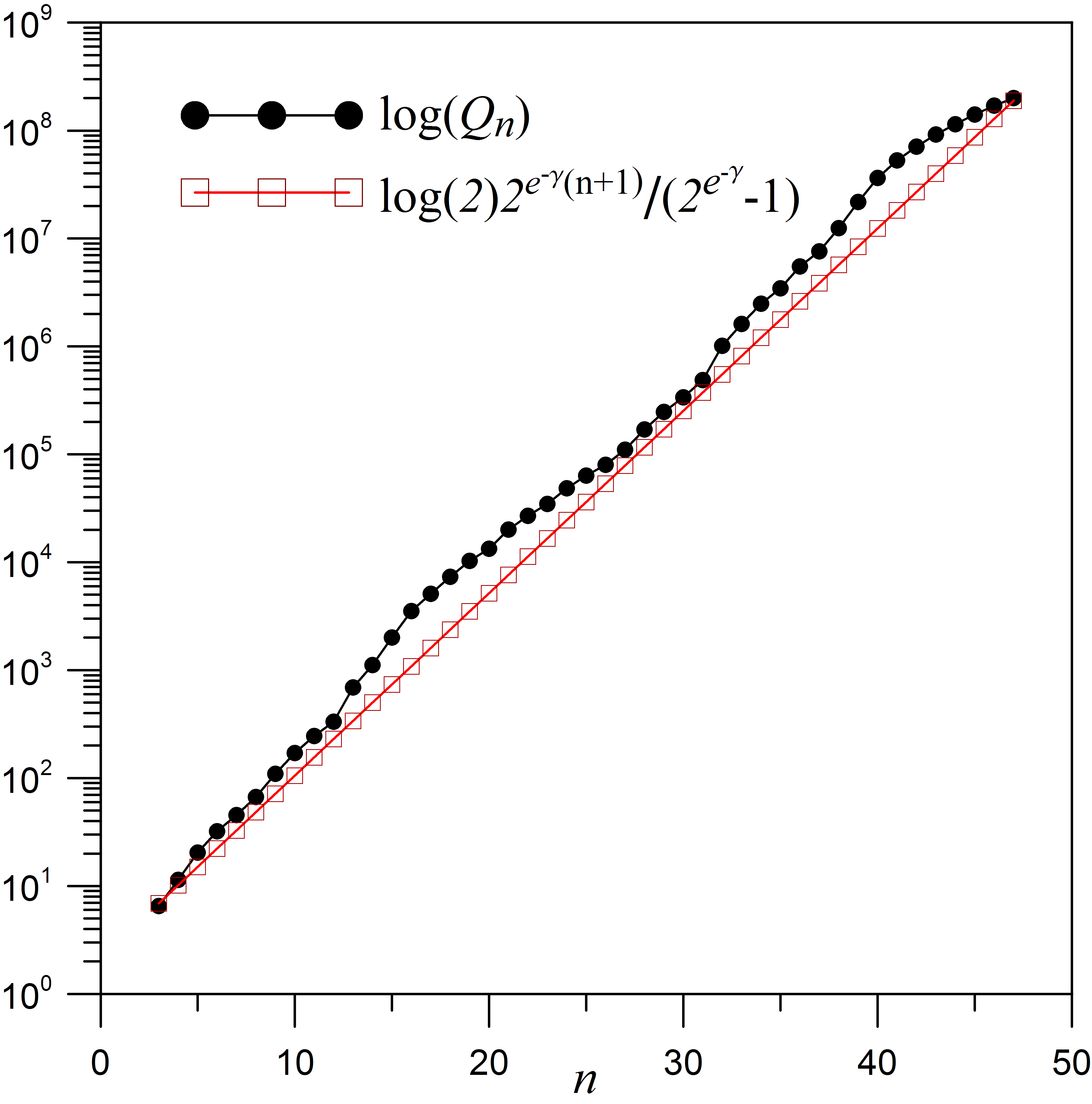} \\
\vspace{0.2cm}
\caption{ Illustration of the inequality (\ref{nierownosc})
for $3\leq n \leq 47$. Although the last points seem to coincide in fact
$Q_{47}= 2.32928\ldots \times 10^{86789810}$, while $2^{c2^{48e^{-\gamma}}}=
1.21513\ldots \times 10^{82034318}$ --- hundreds thousands orders of difference!}
\label{fig-nierownosc}
\end{center}
\end{figure}

Assuming the Wagstaff conjecture $\mathcal{M}_n \sim 2^{2^{n e^{-\gamma}}}$ mentioned
in  the Introduction we obtain that for large $n$
\bee
Q_n > 2^{c2^{(n+1)e^{-\gamma}}},~~~~~~~~c=\frac{1}{2^{e^{-\gamma}}-1}=2.101893933\ldots
\label{nierownosc}
\eee
and thus the transcendence of $u_\mathcal{M}$ will follow from the Davenport--Roth
Theorem (\ref{Davenport-Roth}):
\bee
\frac{\sqrt{\log(n)} \log (\log (Q_n(r)))}{n} \sim \sqrt{\log(n)}\to \infty.
\eee
We illustrate the inequality (\ref{nierownosc})
in Figure \ref{fig-nierownosc} --- the values of labels on the $y$--axis give an idea of the order of
the fast grow of $Q_n(u_{\mathcal{M}})$: the largest for $n=47$ is of the order
$Q_{47}= e^{1.9984\ldots \times 10^8}=2.32928\ldots \times 10^{86789810}$,
see also Table II.

\begin{center}
{\sf TABLE {\bf II}}\\
A sample of values of inverses of the squares of the
$n$-th convergents of $u_\mathcal{M}$ giving an idea of
the speed of convergence of $[0; \mathcal{M}_1, \mathcal{M}_2, \ldots, \mathcal{M}_n]$
to $u_\mathcal{M}$.
\\ \bigskip
\begin{tabular}{|c|c|} \hline
$ n $ & $1/{Q_n^2}$  \\  \hline
 3     & $ 2.131173743\times 10^{-6}  $ \\ \hline
 4     & $   1.320662319\times 10^{-10} $ \\   \hline
 5     & $   1.968416969\times 10^{-18} $ \\   \hline
 6     & $   1.145786956\times 10^{-28} $ \\   \hline
 7     & $   4.168364565\times 10^{-40 } $ \\   \hline
 8     & $   9.038699842\times 10^{-59} $ \\   \hline
 9     & $   1.699990496\times 10^{-95 } $ \\   \hline
$ \vdots $   & $ \vdots $   \\   \hline
 17     & $   9.32543401 \times 10^{-4439 } $ \\   \hline
 18     & $   1.38891910 \times 10^{-6375 } $ \\   \hline
 19     & $   3.81534516 \times 10^{-8936 } $ \\   \hline
 20     & $   4.67942175 \times 10^{-11599 } $ \\   \hline
 $ \vdots $ &  $\vdots $  \\   \hline
 40     & $   4.50116310 \times 10^{-31553835 } $ \\   \hline
 41     & $   5.02100786 \times 10^{-46025300 } $ \\   \hline
 42     & $   3.36434042 \times 10^{-61657758 } $ \\   \hline
 43     & $   3.38166968 \times 10^{-79961861 } $ \\   \hline
 44     & $   2.17906011 \times 10^{-99578575 } $ \\   \hline
 45     & $   5.32688381 \times 10^{-121949118 } $ \\   \hline
 46     & $   1.84595823 \times 10^{-147623244 } $ \\   \hline
 47     & $   1.84313029 \times 10^{-173579621 } $ \\   \hline
\end{tabular}\\
\end{center}

One of the transcendence criterion is the Thue–-Siegel-–Roth Theorem, which we
recall here in the  following form:

{\bf Thue–-Siegel-–Roth Theorem}:  If there exist such $\epsilon>0$ that for
{\it infinitely} many fractions $A_n/B_n$   the inequality 
\bee
%\left|r-\frac{A_n}{B_n}\right|< \frac{1}{B_n^\delta},~~\delta\equiv{2+\epsilon},~~~~n = 1,2,3,...,
\left|r-\frac{A_n}{B_n}\right|< \frac{1}{B_n^{2+\epsilon}},~~~~~n = 1,2,3,...,
\label{Thue–-Siegel-–Roth}
\eee
holds, then $r$ is transcendental.

Let us stress, that $\epsilon$ here does not
depend on $n$  ---   it has to be the same for all fractions $A_n/B_n$.
We can test  the criterion (\ref{Thue–-Siegel-–Roth})  for $u_{\mathcal{M}}$
using as the rational approximations  $A_n/B_n$ the convergents $P_n/Q_n$ of the continued
fraction (\ref{def-uM}).

\begin{figure}
\vspace{-2.3cm}
\hspace{-3.5cm}
\begin{center}
\includegraphics[height=0.4\textheight,angle=0]{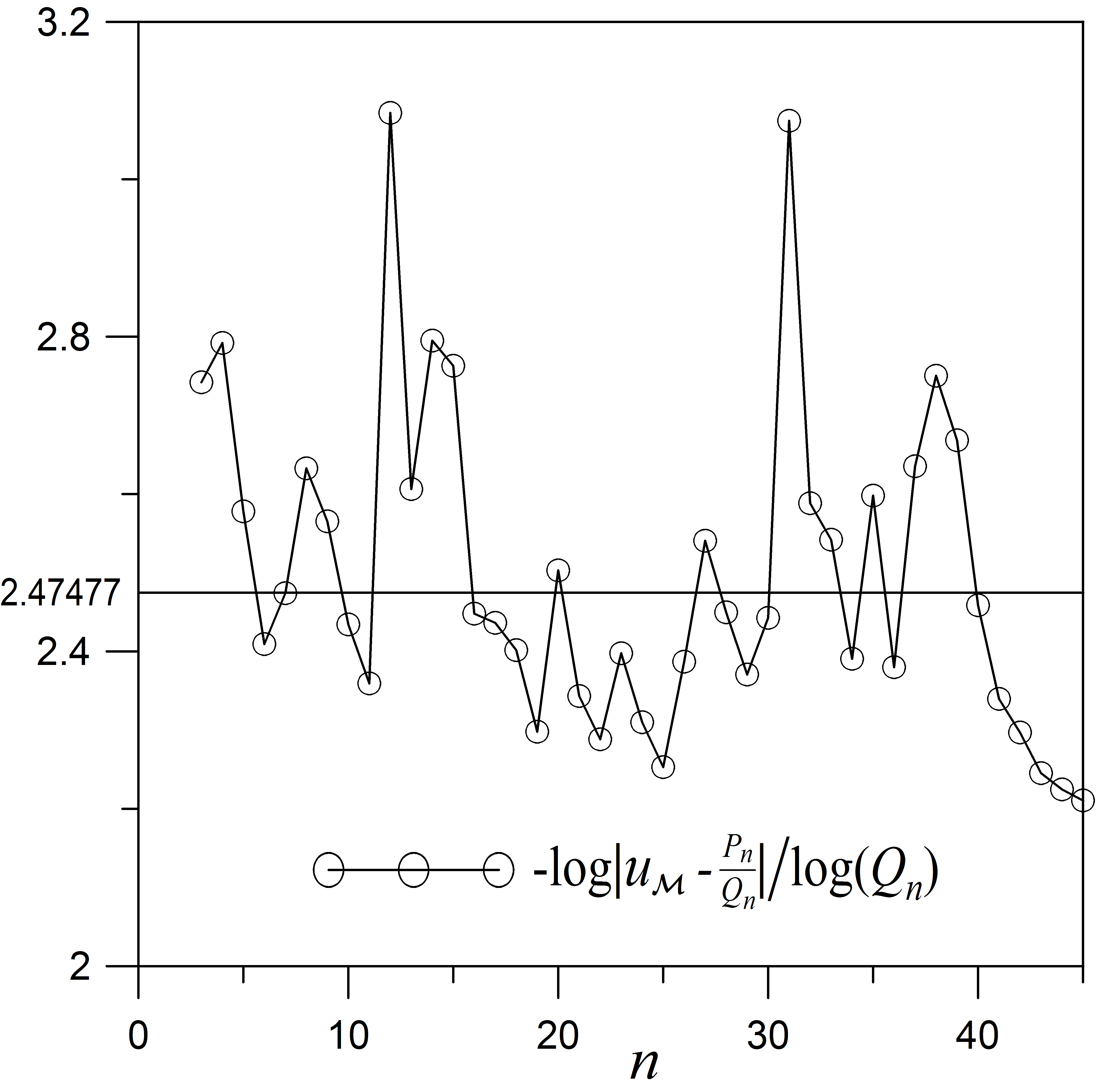} \\
\vspace{0.2cm}
\caption{  The plot of $-\log|u_\mathcal{M}-P_n/Q_n|/\log(Q_n)$
fluctuating around the estimation   (\ref{delta_u_M}). }
\label{fig-delty}
\end{center}
\end{figure}

In \cite{Sondow-2004} J. Sondow has given the estimation  for $\epsilon$
appearing in r.h.s.of (\ref{Thue–-Siegel-–Roth}); namely he proved that for irrational 
numbers with continued fraction expansion $[a_0; a_1, a_2, \ldots]$ and 
convergents $P_n/Q_n$:
\bee
%\delta\leq 1+\limsup_{n \to \infty}\frac{\log Q_{n+1}}{\log Q_n}=2+\limsup_{n \to \infty}\frac{\log a_{n+1}}{\log Q_n}
\epsilon\leq \limsup_{n \to \infty}\frac{\log a_{n+1}}{\log Q_n}.
\eee
Let us denote $\delta\equiv{2+\epsilon}$. From the Wagstaff conjecture we obtain 
that the exponent $\delta$ of $B^\delta$ appearing in on the r.h.s. of 
(\ref{Thue–-Siegel-–Roth}) should be of the order
\bee
\delta \approx 2+2^{e^{-\gamma}}-1=2.47477\ldots  ~~~~~~~~({\rm i.e.} ~~\epsilon\approx 0.47477\ldots)
\label{delta_u_M}
\eee
implying transcendence of $u_\mathcal{M}$. In  Fig.\ref{fig-delty} we present
actual values of  $\delta(u_\mathcal{M}; n)=-\log|u_\mathcal{M}-P_n/Q_n|/\log(Q_n)$
for $n=3, 4, \ldots, 45$ and indeed the values oscillate around 
$1+2^{e^{-\gamma}}=2.47477\ldots$.  First we have calculated $u_\mathcal{M}$ using
all 47 Mersenne primes with accuracy 140000000 digits and  for $n=3, 4, \ldots, 45$
we have calculated convergents $P_n/Q_n$ and next the differences $|u_\mathcal{M} - P_n/Q_n|$
with accuracy $1/Q_n^2$ (see Table II), from which we determined the 
$\delta(u_\mathcal{M}; n)$.  The arithmetic
average of 43 values  $\delta(u_\mathcal{M}; n)$ is $ 2.5002\ldots$, quite
close to the estimated value  (\ref{delta_u_M}). It took a few months of CPU time
to collect data presented in  Fig. \ref{fig-delty}:
It took 12 days of CPU time on
the AMD Opteron 2700 MHz processor to collect data for $n\leq 40$; the  point $n=40$
needed precision of almost 40,000,000 digits, as $|u_\mathcal{M}-P_{40}/Q_{40}|
=1.5033\times 10^{-38789567}$, while $1/Q_{40}^2=4.501\ldots \times 10^{-31553835}$.
To calculate the difference $|u_\mathcal{M}-P_n/Q_n|$ for $n=41, 42, 43$ the precision
of  100000000 digits was needed. For $n=44$ and $n=45$ the difference
$|u_\mathcal{M}-P_n/Q_n|$ was calculated with the precision  130000000 digits 
(see Table II for $n=44$ and $n=45$) and it took about one month of CPU 
time for each  point.

\end{document}